%% file: main.tex
\documentclass[hidelinks,journal]{ieeeconf}

\IEEEoverridecommandlockouts                              %

\usepackage{xr}

\newcommand{\omitted}[1]{}%

\input{titleAndAuthors}

\usepackage{cite}

\usepackage{comment}
\usepackage{siunitx}
\usepackage{relsize}
\usepackage{ifthen}
\usepackage[colorinlistoftodos]{todonotes}

\usepackage[vlined,ruled,linesnumbered]{algorithm2e}
\usepackage{graphics} %
\usepackage{rotating}
\usepackage{color}
\usepackage{enumerate}
\usepackage[T1]{fontenc}
\usepackage{psfrag}
\usepackage{epsfig} %
\usepackage{booktabs}
\usepackage{graphicx,url}
\usepackage{multirow}
\usepackage{array}
\usepackage{latexsym}
\usepackage{amsfonts}
\usepackage{amsmath}
\usepackage{amssymb}
\usepackage{xstring}
\usepackage{algorithmic}
\usepackage{multirow}
\usepackage{xcolor}
\usepackage{prettyref}
\usepackage{flexisym}
\usepackage{bigdelim}
\usepackage{breqn} %
\usepackage{listings}
\let\labelindent\relax
\usepackage{xspace}
\usepackage{bm}
\graphicspath{{./figures/}}
\usepackage{tikz}
\usetikzlibrary{matrix,calc}
\usepackage{lipsum}
\usepackage{mdwlist}

\let\itemize\undefined
\makecompactlist{itemize}{stditemize}
\usepackage{enumitem}
\usepackage{subcaption}

\usepackage{amsthm}
\usepackage{mathtools}

\makeatletter
\let\NAT@parse\undefined
\makeatother
\usepackage{hyperref}
\usepackage{cleveref}

\hypersetup{%
  colorlinks=true,
  linkcolor=blue,
  filecolor=magenta,      
  urlcolor=black,
  citecolor=red,
  linkbordercolor={0 0 1}
}

\input{preamble_symbols.tex}

\input{shortcuts.tex}

\PassOptionsToPackage{end}{algorithmic}

\begin{document}

\maketitle

\thispagestyle{empty}
\pagestyle{empty}

\input{Abstract.tex}

\begin{tikzpicture}[overlay, remember picture]
\path (current page.north east) ++(-4.3,-0.2) node[below left] {
	Please read and cite the journal version of this paper:
    	Z.~Xu, S.S.~Garimella, V.~Tzoumas, 
};
\end{tikzpicture}

\begin{tikzpicture}[overlay, remember picture]
\path (current page.north east) ++(-3.1,-0.6) node[below left] {
    ``Communication- and computation-efficient distributed submodular optimization in robot mesh networks,"
};
\end{tikzpicture}

\begin{tikzpicture}[overlay, remember picture]
\path (current page.north east) ++(-7.7,-1.0) node[below left] {
	IEEE Transactions on Robotics (TRO), 2025.
};
\end{tikzpicture}

\input{1-Introduction.tex}

\input{2-Problem.tex} %

\input{3-Algorithms.tex} %

\input{4-Resources.tex} %

\input{5-Guarantees.tex} %

\input{6-Experiments.tex}

\input{7-Conclusion.tex}

\section*{Acknowledgements}
We thank Robey et al.~\cite{robey2021optimal} and Konda et al.~\cite{konda2021execution} for sharing with us the code of their numerical evaluations. We also thank Hongyu Zhou of the University of Michigan for providing comments on the paper.

\bibliographystyle{IEEEtran}
\bibliography{references,security_references}

\end{document}

%% file: titleAndAuthors.tex
\title{\fontsize{22}{22} \selectfont Resource-Aware Distributed Submodular Maximization: \\A Paradigm for Multi-Robot Decision-Making
}
\author{Zirui Xu,$^\star$ Vasileios Tzoumas$^\dagger$
\thanks{Submitted for publication in March 2022 and accepted to CDC~2022.}
	\thanks{$^\star$Department of Aerospace Engineering, University of Michigan, Ann Arbor, MI 48109 USA;  \href{mailto:ziruixu@umich.edu}{\tt\footnotesize ziruixu@umich.edu}}
	\thanks{$^{\dagger}$Department of Aerospace Engineering and Robotics Institute, University of Michigan, Ann Arbor, MI 48109 USA; \href{mailto:vtzoumas@umich.edu}{\tt\footnotesize vtzoumas@umich.edu}}
}

%% file: preamble_symbols.tex
\newtheorem{theorem}{Theorem}
\newtheorem{problem}{Problem}

\newtheorem{assumption}{Assumption}
\newtheorem{definition}{Definition}
\newtheorem{proposition}{Proposition}

\newtheorem{remark}{Remark}

\newcommand{\bdmath}{\begin{dmath}}
\newcommand{\edmath}{\end{dmath}}
\newcommand{\beq}{\begin{equation}}
\newcommand{\eeq}{\end{equation}}
\newcommand{\bdm}{\begin{displaymath}}
\newcommand{\edm}{\end{displaymath}}
\newcommand{\bea}{\begin{eqnarray}}
\newcommand{\eea}{\end{eqnarray}}
\newcommand{\beal}{\beq \begin{array}{lll}}
\newcommand{\eeal}{\end{array} \eeq}
\newcommand{\beas}{\begin{eqnarray*}}
\newcommand{\eeas}{\end{eqnarray*}}
\newcommand{\ba}{\begin{array}}
\newcommand{\ea}{\end{array}}
\newcommand{\bit}{\begin{itemize}}
\newcommand{\eit}{\end{itemize}}
\newcommand{\ben}{\begin{enumerate}}
\newcommand{\een}{\end{enumerate}}

\newcommand{\calA}{{\cal A}}
\newcommand{\calB}{{\cal B}}
\newcommand{\calC}{{\cal C}}

\newcommand{\calE}{{\cal E}}

\newcommand{\calG}{{\cal G}}

\newcommand{\calI}{{\cal I}}

\newcommand{\calN}{{\cal N}}

\newcommand{\calS}{{\cal S}}

\newcommand{\calV}{{\cal V}}

\newcommand{\calX}{{\cal X}}

\definecolor{myblue}{RGB}{65 105 225}

\newcommand{\hide}[1]{}

\newcommand{\hiddenText}{{\color{gray} hidden text.}}
\newcommand{\hideWithText}[1]{\hiddenText}

\newcommand{\opt}{^{\star}}

\newcommand{\scenario}[1]{{\smaller \sf#1}\xspace}

%% file: shortcuts.tex
\newcommand{\ie}{\emph{i.e.},\xspace}
\newcommand{\eg}{\emph{e.g.},\xspace}
\newcommand{\myin}{\, \in \,}

\newcommand{\alg}{\scenario{RAG}}

\newcommand{\myprob}{Problem~\ref{pr:main}\xspace}

\newcommand{\myParagraph}[1]{{\bf #1.}\xspace}

\newcommand{\inneigh}{\calN^-}
\newcommand{\outneigh}{\calN^+}
\renewcommand{\opt}{\scenario{OPT}}
\newcommand{\single}{s}

\newcommand{\ourcurv}{\scenario{coin}}

\newcommand{\ouragent}{\mathcal{N}}
\newcommand{\oursol}{\calS^{\alg}}
\newcommand{\singlesol}{s^{\alg}}

\newcommand{\numl}{\scenario{length}_\#}
\newcommand{\el}{\scenario{length}_s}
\newcommand{\rob}{j}

\newcommand{\agentselected}{\calI}
\newcommand{\solutionselected}{\calS}

%% file: Abstract.tex
\begin{abstract}
We introduce {the first} algorithm for distributed decision-making that provably balances the trade-off of \textit{centralization}, for {global near-optimality}, vs.~\textit{decentralization}, for near-minimal on-board computation, communication, and memory resources.  
We are motivated by the future of autonomy that involves \textit{heterogeneous} robots collaborating in complex~tasks, such as image covering, target tracking, and area monitoring.  Current algorithms, such as consensus algorithms, are \textit{insufficient} to fulfill this future: they achieve distributed communication \textit{only}, at the expense of high communication, computation, and memory overloads. 
{A shift to \textit{resource-aware} algorithms is needed, that can account for each robot's on-board resources, \textit{independently}.}
We provide the first resource-aware algorithm, \textit{Resource-Aware distributed Greedy} (\alg). We focus on maximization problems involving monotone and 2nd-order submodular functions, a diminishing returns property.
\alg has {near-minimal} on-board resource requirements. Each agent can afford to run the algorithm by adjusting the size of its neighborhood, even if that means selecting actions in complete isolation. 
\alg has {provable} approximation performance, where each agent can {independently} determine its contribution.
All in all, \alg is the first algorithm to quantify the trade-off of \emph{centralization}, for global near-optimality, vs.~\emph{decentralization}, for near-minimal on-board  resource requirements.
To capture the trade-off, we introduce the notion of \emph{Centralization Of Information among non-Neighbors} (\scenario{COIN}).  We validate \alg in simulated scenarios of image covering with mobile robots.
\end{abstract}

%% file: 1-Introduction.tex
\vspace{-10mm}
\section{Introduction}\label{sec:Intro}
In the future, robots with heterogeneous on-board capabi- lities will be teaming up
to complete complex tasks such as:
\begin{itemize}
    \item \textit{Image Covering}: How swarms of tiny to large robots can collaboratively map (image cover) an unknown environment, such as an earthquake-hit building?~\cite{mcguire2019minimal}
    \item \textit{Area Monitoring}: How swarms of ground and air robots can collaboratively monitor an environment to detect rare events, such as fires in a forest?~{\cite{kumar2004robot}}
    \item \textit{Target Tracking}: How a large-scale distributed network of air and space vehicles can coordinate their motions to track multiple evading targets over a large area?~{\cite{corah2021scalable}}
\end{itemize}
The robots' heterogeneous capabilities (speed, size, on-board cameras, etc.) offer tremendous advantages in all aforementioned tasks: {for example}, in the image covering scenario, the tiny robots offer the advantage of agility, being able to navigate narrow spaces in earthquake-hit buildings; and the larger robots offer the advantages of reliability, being able to carry larger and higher-resolution cameras, for longer.

\definecolor{lightgray}{rgb}{0.83, 0.83, 0.83}
{\footnotesize
\begin{table}[t]
\captionsetup{font=small}
\centering
\renewcommand{\arraystretch}{1.3}
\begin{tabular}[t]{ | l || l |  }
 \multicolumn{2}{c}{}\\
\hline
\hline
 \multicolumn{2}{c}{\textbf{
  \fontsize{9}{9} \selectfont 
  \!\!\hspace{-2mm}Resource-Aware Distributed Decision-Making: A Paradigm\!\!}} \\ 
 \hline
 \hline
 \textbf{Computations per Agent} &  \begin{tabular}{@{}l@{}} 
 Proportional to the size of the agent's \vspace{-0.75mm}\\ available action set \end{tabular} \\ 
 \hline 
 \textbf{Communication Rounds} & Proportional to the number of agents \\   
 \hline 
 \textbf{Memory per Message} &  Length of a real number or an action \\
 \hline
 \textbf{Communication Topology} & Directed and even disconnected \\
 \hline
 \textbf{Suboptimality Guarantee} &  \begin{tabular}{@{}l@{}}
 Gracefully balances the trade-off of \vspace{-0.75mm}\\centralization vs.~decentralization \end{tabular} \\
 \hline
\end{tabular}
\renewcommand{\arraystretch}{1}
\caption{\textbf{Resource-Aware Distributed Optimization Paradigm.} We define \emph{resource-awareness} in distributed optimization across five performance pillars.  The pillars define an algorithm that (i\hspace{1pt}---\hspace{1pt}see first 3 pillars) has  \emph{minimal} computation, communication and memory-storage requirements, and is {affordable} by any agent when the agent chooses a small enough neighborhood to meet its 
resources; (ii\hspace{1pt}---\hspace{1pt}see the 4th pillar) is applicable to even disconnected communication topologies, which is required when agents lack or are diminished of on-board resources for information exchange; and (iii\hspace{1pt}---\hspace{1pt}see the 5th pillar) has a suboptimality guarantee that balances 
the trade-off of \emph{centralization}, for global near-optimality, vs.~\emph{decentralization}, for near-minimal on-board resource requirements.
} \label{tab:intro}
\end{table}
}

But heterogeneity in capabilities also implies heterogeneity in on-board resources: {for example,} tiny robots have limited computation, communication, and memory resources~\cite{mcguire2019minimal}.  Thus, {mere} distributed \textit{communication} among the robots is \textit{insufficient} for the success of their tasks.  Instead, a holistic, \textit{resource-aware} distributed collaboration is necessitated
that respects each robot's on-board capacity for computation, communication, and memory storage. 

Current algorithms, such as consensus algorithms, are \emph{insufficient} to meet the need for resource-awareness: they achieve distributed communication but at the expense of communication, computation, and memory overloads.
Hence, robots with limited on-board resources, such as the tiny (27 grams) Crazyflies, cannot afford to use the algorithms~\cite{mcguire2019minimal}.  {Also,} real-time performance is compromised by the latency caused by the computation and communication~overloads.

In this paper, we shift focus from the current \textit{distributed-communication} only optimization paradigm to the \emph{resource-aware} distributed optimization paradigm in \Cref{tab:intro}.   

We focus on scenarios where the robots' tasks are captured by an objective function
that is monotone and 2nd-order submodular~{\cite{crama1989characterization,foldes2005submodularity}}, a diminishing returns property.   Such functions appear  in tasks of image covering {\cite{corah2018distributed} and vehicle deployment \cite{downie2022submodular}}, among others. Then, the aforementioned  tasks require the robots to distributively 
\begin{equation}\label{eq:intro}
\max_{\single_i\myin\mathcal{V}_i, \, \forall\, i\myin \calN} \; f(\,\{s_i\}_{i\myin \calN}\,),
\end{equation}
where $\calN$ is set of agents/robots, $s_i$ is agent $i$'s action, $\calV_i$ is agent $i$'s set of available actions (\eg motion primitives), and $f:2^{\prod_{i \in \calN}\,\calV_i}\mapsto\mathbb{R}$ is the objective function (\eg total area covered by agents' cameras at the current time step).  In online settings, the robots may need to solve a new version of \cref{eq:intro} at each time step \mbox{(\eg in a receding-horizon fashion).} 

\myParagraph{Related Work}  
Problem~\ref{eq:intro} is NP-hard, being combinatorial, even when $f$ is monotone and  submodular~{\cite{Feige:1998:TLN:285055.285059}}.  It has been actively researched in the last 40 years
in the optimization, control, and operations research literature~\cite{fisher1978analysis,crama1989characterization,Feige:1998:TLN:285055.285059, foldes2005submodularity,krause08efficient,calinescu2011maximizing,wang2015accelerated,atanasov2015decentralized,mirzasoleiman2016distributed,ramalingam2017efficient,roberts2017submodular,sviridenko2017optimal,gharesifard2017distributed,grimsman2018impact,corah2018distributed,corah2019distributed,du2020jacobi,robey2021optimal,rezazadeh2021distributed,mokhtari2018decentralized,konda2021execution,downie2022submodular,chen2022higher,corah2021scalable}.  Although near-optimal approximation algorithms have been achieved,
they focus on distributed communication only, instead of a holistically resource-aware optimization per \Cref{tab:intro}.  For example, the 
\textit{sequential greedy}~\cite{fisher1978analysis} and its variants~{\cite{sviridenko2017optimal,gharesifard2017distributed,grimsman2018impact,corah2018distributed,konda2021execution}}, 
require increasing memory storage (agents act sequentially, and each agent passes the information about \textit{all} previous agents to the next). 
Further the \textit{consensus-based} distributed algorithms~{\cite{mokhtari2018decentralized,robey2021optimal,rezazadeh2021distributed}}, 
although they achieve distributed communication,
require excessive computations and communication rounds per agent.  
No current algorithm has suboptimality guarantees for even disconnected communication networks, nor captures 
the trade-off of \emph{centralization}, for global near-optimality, vs.~\emph{decentralization}, for near-minimal on-board resource requirements.

\myParagraph{Contributions} 
We shift focus to holistically resource-aware distributed optimization.  Along with novel definitions and theory, we provide the first algorithm
balancing the  trade-off of \emph{centralization}, 
for global near-optimality,
vs.~\emph{decentralization}, for near-minimal on-board resource requirements.

\myParagraph{1.~\!\!Resource-Aware Optimization Paradigm} Our first contribution is the definition of a \emph{resource-aware} distributed algorithm.  The definition is summarized in
\Cref{tab:intro}.

\myParagraph{2.~\!\!Resource-Aware Algorithm} 
We introduce the first resource-aware distributed algorithm for \cref{eq:intro} (\Cref{sec:algorithm}). 
Per \Cref{tab:intro}: 
(i) \alg has \emph{near-minimal} computation, communication, and memory requirements (\Cref{sec:resource}).
(ii) Each agent can afford to run \alg by adjusting the size of its neighborhood, even if that means selecting actions in complete isolation.
(iii) \alg enjoys a suboptimality bound that captures 
the trade-off of \emph{centralization}, for global near-optimality, vs.~\emph{decentralization}, for near-minimal on-board resource requirements 
(\Cref{sec:guarantee}). All agents can {independently} decide their contribution to the approximation performance by choosing the size of their neighborhood, and accounting at the same time for their on-board resources.
\alg is the first algorithm to quantify the trade-off.

\myParagraph{3.~\mbox{Centralization of Information}} 
We introduce the notion of \emph{Centralization of Information among non-Neighbors} (\scenario{COIN}) to quantify \alg's suboptimality.  \scenario{COIN} captures the information overlap between an agent and its {non-neighbors}.  

\myParagraph{Evaluation on Robotic Application} We evaluate \alg in simulated scenarios of image covering with mobile robots (\Cref{sec:experiments}). We first compare \alg with the state of the art.  Then, we evaluate the trade-off of centralization vs.~decentralization with respect to \alg's performance.  To enable the comparison with the state of art, we assume undirected and connected networks.
\alg demonstrates superior or comparable performance, requiring, \eg (i) $5$ orders of magnitude less computation time vs.~the state-of-the-art consensus algorithm in~\cite{robey2021optimal}, (ii) $1$ order of magnitude fewer communication rounds vs.~the consensus algorithm in~\cite{robey2021optimal}, and comparable communication rounds vs.~the greedy in~\cite{konda2021execution}, and (iii) the least memory (\eg $40\%$ less vs.~the consensus algorithm in~\cite{robey2021optimal}).  Still, \alg has the best approximation performance.

%% file: 2-Problem.tex
\section{Distributed Submodular Maximization: \\A Multi-Robot Decision-Making  Perspective
}\label{sec:problem}

We define the Distributed Submodular Maximization problem of this paper (\Cref{pr:main}).  We use the notation:
\begin{itemize}
    \item  $\calG\triangleq \{\calN, \calE\}$ is a communication network 
    with nodes $\calN$ and edges $\calE$.  Nodes represent \emph{agents} (\eg robots), and edges represent \emph{communication channels}. 
    \item $\calN^-_i\triangleq \{j\in \calN : (j,i) \in \calE\}$, for all $i\in\calN$; \ie $\inneigh_i$ is the in-neighbors of $i$.
    \item $\calN^+_i\triangleq \{j\in \calN : (i,j) \in \calE\}$, for all $i\in\calN$, given graph $\calG$; \ie $\outneigh_i$ is the out-neighbors of $i$.  If $\calG$ is undirected, then $\inneigh_i=\outneigh_i$, for all $i\in\calN$.
    \item $d(i,j)$ is a distance metric between an $i\in \calN$ and a $j \in \calN$; \eg $d(i,j)$ may be the Euclidean distance between $i$ and $j$, when $i$ and $j$ are robots in the $\scenario{3D}$ space. 
    \item $\calV_\calN \triangleq \prod_{i\myin \calN} \,\calV_i$ given a collection of sets $\{\calV_i\}_{i\myin \calN}$; \ie  $\calV_\calN$ is the cross-product of the sets in $\{\calV_i\}_{i\myin \calN}$; 
    \item $f(s\,|\,\calA)\triangleq f(\calA \cup \{s\})-f(\calA)$, given a set function $f:2^\calV\mapsto \mathbb{R}$,  $s \in \calV$, and $\calA \subseteq\calV$; \ie $f(s\,|\,\calA)$ is the marginal gain in $f$ for adding $s$ to $\calA$.
\end{itemize}

The following preliminary framework is also required.

\myParagraph{Agents} The terms ``\emph{agent}'' and ``\emph{robot}'' are used interchangeably in this paper.  $\calN$ is the set of all robots.  The robots cooperate towards a task, such as \emph{image covering}. 
$\calV_i$ is a \textit{discrete} set of actions available to each robot $i$.  For example, in image covering, $\calV_i$ may be the set of motion primitives that robot $i$ can execute to move in the environment.

\myParagraph{Communication Network} The communication network $\calG$ among the robots may be directed and even disconnected.  If $(j,i)\in\calE$, then a communication channel exists from robot $j$ to robot $i$: $i$ can receive, store, and process the information from $j$.  The set of all robots that can send information to $i$ is $\calN^-_i$, \ie $i$'s in-neighborhood.   The set of all robots that $i$ can send information to is $\calN^+_i$, \ie $i$'s out-neighborhood. 

\begin{remark}[Resource-aware in-neighborhood selection based on information overlap] \label{rem:info-overlap} In this paper, $\calN^-_i$ has been implicitly decided by robot $i$ given $i$'s on-board resources
and based on a \emph{distance metric} $d(i,j)$ capturing the \emph{information overlap} between $i$ and any robot within communication range. Particularly, $d(i,j)$ is considered to \emph{increase} as the {information overlap} \emph{drops}.  
\emph{E.g.}, in \emph{image covering},
since each camera's field of view is finite, the distance metric $d$ can be proportional to the physical distance of the robots.  When $d(i,j)$ is sufficiently large, such that the field of view of robots $i$ and $j$ stop overlapping, then the robots know that their information is non-overlapping, and may stop exchanging information to reserve on-board resources.
\end{remark}

\begin{remark}[Resource-aware out-neighborhood selection] \label{rem:out-neigh} In this paper, $\calN^+_i$ has been implicitly decided by robot $i$ given robot $i$'s on-board resources.  
\end{remark}

\myParagraph{Objective Function} The robots coordinate their actions to maximize an objective function.  In 
tasks, such as image covering, target tracking, and persistent monitoring, typical objective functions are the \emph{covering functions} \cite{corah2018distributed,robey2021optimal, downie2022submodular}.  Intuitively, these functions capture how much area{/information} is covered given the actions of all robots.  
They satisfy the properties defined below (\Cref{def:submodular} and \Cref{def:conditioning}).

\begin{definition}[Normalized and Non-Decreasing Submodular Set Function{~\cite{fisher1978analysis}}]\label{def:submodular}
A set function $f:2^\calV\mapsto \mathbb{R}$ is \emph{normalized and non-decreasing submodular} if and only if 
\begin{itemize}
\item $f(\emptyset)=0$;
\item $f(\calA)\leq f(\calB)$, for any $\calA\subseteq \calB\subseteq \calV$;
\item $f(s\,|\,\calA)\geq f(s\,|\,{\mathcal{B}})$, for any $\calA\subseteq {\mathcal{B}}\subseteq\calV$ and $s\in \calV$.
\end{itemize}
\end{definition}

Normalization $(f(\emptyset)=0)$ holds without loss of generality.  In contrast, monotonicity and submodularity are intrinsic to the function.  
Intuitively, if $f(\calA)$ captures the area \emph{covered} by a set $\calA$ of activated cameras, then the more sensors are activated $(\calA\subseteq \calB)$, the more area is covered $(f(\calA)\leq f(\calB))$; this is the non-decreasing property.  Also, the marginal gain of covered area caused by activating a camera $s$ \emph{drops} $(f(s\,|\,\calA)\geq f(s\,|\,{\mathcal{B}}))$ when \emph{more} cameras are already activated $(\calA\subseteq \calB)$; this is the submodularity~property.

\begin{definition}[2nd-order Submodular Set Function{~\cite{crama1989characterization,foldes2005submodularity}}]\label{def:conditioning}
$f:2^\calV\mapsto \mathbb{R}$ is \emph{2nd-order submodular} if and only if 
\begin{equation}\label{eq:conditioning}
    f(s\,|\,\calC) - f(s\,|\,\calA\cup\calC) \geq f(s\,|\,\calB\cup\calC) - f(s\,|\,\calA\cup\calB\cup\calC),
\end{equation}
for any \emph{disjoint} $\calA, \calB, \calC\subseteq \calV$ ($\calA \cap \calB \cap \calC =\emptyset$) and $s\in\calV$.
\end{definition}

The 2nd-order submodularity is another intrinsic property to the function.  
Intuitively, if $f(\calA)$ captures the area \emph{covered} by a set $\calA$ of cameras, then \emph{marginal gain of the marginal gains} drops when more cameras are already activated.

\myParagraph{Problem Definition}  In this paper, we focus on:

\begin{problem}[Distributed Submodular Maximization]\label{pr:main}
Each robot $i \in \calN$ independently selects an action $s_i$, \emph{upon receiving information from and about the in-neighbors $\calN^-_i$ only,}
such that the robots' actions $\{s_i\}_{i \myin \calN}$ solve the 
\begin{equation}\label{eq:problem}
\max_{\single_i\myin\mathcal{V}_i, \, \forall\, i\myin \calN} \; f(\,\{s_i\}_{i\myin \calN}\,),
\end{equation}
where $f:2^{\calV_\calN}\mapsto \mathbb{R}$ is a normalized, non-decreasing submodular, and 2nd-order submodular set function. 
\end{problem}

\myprob requires each robot $i$ to independently choose an action $s_i$ to maximize the global objective $f$, only based on \emph{local communication \emph{and} information}.  Particularly, if after $\tau_i$ communication rounds (\ie iterations of information exchange) robot $i$ decides to select an action $s_i$, then {$$s_i = \phi_i\!\left(\!\begin{tabular}{@{}l@{}}
	${\color{white}\{}$Information received from robot $i$'s sensors,\\
    $\{$Information received from robot $j$  till $\tau_i$$\}_{j\myin \calN^-_i}$
	\end{tabular}\!\right)$$}for a decision algorithm $\phi_i$ to be found in this paper. 
	
\begin{assumption}\label{ass:info}
Each robot $i$ may receive the action $s_j$ of an in-neighbor $j \in \calN^-_i$ as information, and, then, robot $i$ can locally compute the marginal gain of any of its own $s_i \in \calV_i$. 
\end{assumption}

The assumption is common in all distributed optimization algorithms in the literature {\cite{gharesifard2017distributed,grimsman2018impact,mokhtari2018decentralized,robey2021optimal,downie2022submodular,du2020jacobi,rezazadeh2021distributed,corah2018distributed,konda2021execution,liu2021distributed,corah2019distributed,corah2021scalable}}. 
In practice, robot $j$ needs to transmit to robot $i$, along with $j$'s action, all other required information such that $i$ can compute locally the marginal gains of any of its own candidate actions. 

In Section~\ref{sec:algorithm} we introduce \alg, an algorithm that runs locally on each robot $i$, playing the role of $\phi_i$. %

\begin{assumption}\label{ass:fixed}
The communication network $\calG$ is fixed between the executions of the algorithm.
\end{assumption}

That is, we assume no communication failures once the algorithm has started, and till its end.  Still, $\calG$ can be \emph{dynamic} across consecutive time-steps $t=1, 2, \ldots$\;, when \myprob is applied in a \emph{receding-horizon fashion} (Remark~\ref{rem:rhc}). 

\begin{remark}[Receding-Horizon Control, and Need for Minimal Communication and Computation] \label{rem:rhc}
Image covering, target tracking, and persistent monitoring are dynamic tasks that require the robots to react across consecutive time-steps $t=1, 2, \ldots$\;.  Then, \myprob must be solved in a \emph{receding-horizon fashion}~{\cite{camacho2013model}}. This becomes possible only if the time interval between any two consecutive steps $t$ and $t+1$ 
can contain 
the required number of communication rounds to solve \Cref{pr:main}.  Thus, for faster reaction, the smaller the number of communication rounds must be, and the 
smaller the computation effort per round must be.  
Otherwise, real-time performance will be compromised by latency. 
\end{remark}

%% file: 3-Algorithms.tex
\section{Resource-Aware distributed Greedy (\alg) Algorithm} \label{sec:algorithm}

{
\setlength{\textfloatsep}{3mm}
\begin{algorithm}[t!]
	\caption{\mbox{Resource-Aware distributed Greedy (\alg).}
	}
	\begin{algorithmic}[1]
		\REQUIRE\!Agent $i$'s action set\hspace{-1.25pt} $\mathcal{V}_i$; in-neighbors\hspace{-1.25pt} set $\calN_i^-$; out-neighbors\hspace{-1.25pt} set $\calN_i^+$; \hspace{-2pt}normalized, \hspace{-1.25pt}non-decreasing \hspace{-1.25pt}submodular, and 2nd-order submodular set function $f:2^{\mathcal{V}_\calN} \mapsto \mathbb{R}$. %
		\ENSURE Agent $i$'s action $\singlesol_i$.%
		\medskip
		
		\STATE $\agentselected_i\leftarrow\emptyset$;~~~$\solutionselected_i\leftarrow\emptyset$;~~~$\singlesol_i\leftarrow\emptyset$;\label{line:initiliaze} {\color{gray}// \!$\calI_i$ \!stores the agents\\ in $\calN^-_i$ that have selected an action; $\calS_i$ stores $\calI_i$'s selected actions; $s_i^\alg$ stores agent $i$'s selected action}
        \WHILE {$\singlesol_i = \emptyset$}
        \STATE $s_i\gets \arg\max_{y \myin \mathcal{V}_i}\,f(\solutionselected_i\cup \{y\})-f(\solutionselected_i)$;
        \STATE $g_i\gets f(\solutionselected_i\cup \{s_i\})-f(\solutionselected_i)$;
        \STATE \textbf{transmit} $g_i$ \textbf{to} each agent $j\in\mathcal{N}_i^+\setminus\agentselected_i$ and \textbf{receive} $\{g_\rob\}_{\rob\myin{\mathcal{N}_i^-\setminus\agentselected_i}}$;
        \IF {$i=\arg\max_{\rob \myin \mathcal{N}_i^-\cup\{i\}\setminus\agentselected_i}\,g_{\rob}$}
        \STATE $\singlesol_i \leftarrow s_i$; {\color{gray}// $i$ selects  action}
        \STATE \textbf{transmit} $s_i$ \textbf{to} each agent $j \in \mathcal{N}_i^+\setminus\agentselected_i$; {\color{gray}// $i$ has \\the  best action across $\{\mathcal{N}_i^-\cup\{i\}\}\setminus\agentselected_i$}
        \ELSE
        \STATE \textbf{denote} by $\agentselected^{\scenario{new}}_i$ the set of agent(s) $j \in \calN_i^-\setminus \calI_i$ \\ that selected action(s) in this iteration;
        \STATE \textbf{receive} $\singlesol_j$ 
        \textbf{from} each agent $j \in \agentselected^{\scenario{new}}_i$;
        \STATE $\agentselected_i \leftarrow \agentselected_i \cup \agentselected^{\scenario{new}}_i$;
        \STATE $\solutionselected_i \leftarrow\solutionselected_i \cup \{\singlesol_j\}_{j\myin\calI_i^{\mathlarger{\scenario{new}}}}$;
        \ENDIF
		\ENDWHILE\label{line:end_for_3}
        \RETURN $\singlesol_i$.
	
	\end{algorithmic}\label{alg:dec_sub_max}
\end{algorithm}
}

We present the \textit{Resource-Aware distributed Greedy} (\alg) algorithm, the first resource-aware algorithm for  \Cref{pr:main}.  

\alg's pseudo-code is given in \Cref{alg:dec_sub_max}.  The algorithm requires only \emph{local} information exchange, \emph{among} only neighboring robots and \emph{about} only neighboring robots.  Each agent $i$ starts by selecting the action $s_i$ with the largest marginal gain $g_i$ (lines 3--4).
Then, instead of exchanging information with all other agents $\ouragent\setminus\{i\}$, $i$ exchanges information only with the in-neighbors $\calN_i^-$ and out-neighbors $\calN_i^+$ (line 5).  Afterwards, $i$ checks whether $i=\arg\max_{\rob \myin \mathcal{N}_i^-\cup\{i\}}g_{\rob}$, \ie whether $i$ is the ``best agent'' among the in-neighbors only, instead of all agents in the network 
(line 6). If yes, then $i$ selects $s_i$ as its action (line 7), and lets only its out-neighbors know its selection (line 8). Otherwise, $i$ receives the action{(}s{)} from the agent{(}s{)} $j \in \calN_i^-\setminus \calI_i$ that just selected action{(}s{)} in this iteration ({\ie} set $\calI_i^{\scenario{new}}$ in line 10), and continues onto the next iteration (lines 9--15). Notably, $\calI_i^{\scenario{new}}$ may contain multiple agents, and may even be empty.

\begin{remark}[Directed and Disconnected Communication Topology]\label{rem:disconnected} \alg is valid for directed and even disconnected communication topologies, in accordance to the paradigm \Cref{tab:intro}.  For example, if $\calN_i^-=\calN_i^+=\emptyset$ in Algorithm~\ref{alg:dec_sub_max}, then agent $i$ is completely disconnected from the network. 
\end{remark}

%% file: 4-Resources.tex
\section{Computation, Communication, and Memory Requirements of \alg}\label{sec:resource}

\definecolor{OliveGreen}{rgb}{0,0.5,0}
{\footnotesize
\begin{table*}[t]
\captionsetup{font=footnotesize}
\centering
\renewcommand{\arraystretch}{1.5}
\begin{tabular}{ |l||c|c|c||c| }
 \hline
  & \multicolumn{3}{c||}{\textbf{Continuous Domain}} & \textbf{} \\
 \hline
  \textbf{Method} & {Du et al. \cite{du2020jacobi}} & {Robey et al. \cite{robey2021optimal}} & {Rezazadeh and Kia \cite{rezazadeh2021distributed}} & \textbf{\alg} (this paper) \\
 \hline 
 \hline
 \textbf{Computations per Agent} & {$\sim\Theta(M\,|\ouragent|^2)$} & {$\Omega(M\,|\ouragent|^{2.5}\,/\,\epsilon)$} &\!$\Omega(M\,|\ouragent|^2\,\scenario{diam}(\calG)\,/\,\epsilon)$\!& $O(|\mathcal{V}_i|\,|\mathcal{N}^-_i|)$ \\ 
 \hline
 \textbf{Communication Rounds} &
  {$\sim\Theta(|\ouragent|^2)$} & $\Omega(|\ouragent|^{2.5}\,/\,\epsilon)$ & $\Omega(|\ouragent|^2\,\scenario{diam}(\calG)\,/\,\epsilon)$ & $\leq \,2\, |\ouragent|\,-\,2$ \\ 
 \hline
 \textbf{Memory per Message} & $M\,\el$ & {\color{black}$|\mathcal{V}_\calN|\,\numl$} & $|\mathcal{V}_\calN|\,\numl$ & {\color{black}$\max{(\numl,\,\el)}$} \\ 
 \hline
 \textbf{Communication Topology} & connected, undirected & connected, undirected & connected, undirected & {\color{black}even disconnected, directed} \\
 \hline
 \textbf{Suboptimality Guarantee} & $(1/2-\epsilon)\,\opt$ & $(1-1/e)\,\opt$ $-\,\epsilon$ & $(1-1/e-\epsilon)\,\opt$ &\!$1/2\,(\opt-{\sum_{i\in \calN}{\ourcurv}_i})$\!\\
 \hline
  \hline
  \hline
  & \multicolumn{3}{c||}{\textbf{Discrete Domain}} & \textbf{} \\
 \hline
 \textbf{Method} & {Corah and Michael \cite{corah2018distributed}} & {Liu et al. \cite{liu2021distributed}} & {Konda et al. \cite{konda2021execution}} & \textbf{\alg} (this paper) \\
 \hline 
 \hline
 \textbf{Computations per Agent} & $O(|\mathcal{V}_i|\,| {\mathcal{V}_{\mathcal{N}\setminus\{i\}}}|)$ & $O(|\mathcal{V}_i|\,|\ouragent|^2)$ & {\color{black}$|\mathcal{V}_i|$} & $O(|\mathcal{V}_i|\,|\mathcal{N}^-_i|)$ \\
 \hline
 \textbf{Communication Rounds} & {\color{black}$\Omega(1\,/\,\epsilon)$} & $\leq\,(2\,|\ouragent|\,+\,2)\,\scenario{diam}(\mathcal{G})$ & ${\leq}\,2\,|\ouragent|\,-\,2$ & {\color{black}$\leq\,2\,|\ouragent|\,-\,2$} \\
 \hline
 \textbf{Memory per Message} &\!$\max(|\mathcal{V}_i|,\,\Omega(1/\epsilon))\,\el$\!&\!$|\ouragent|(\numl\,+\,\el)$\!& $(|\ouragent|\,-\,1)\,\el$ & {\color{black}$\max{(\numl,\,\el)}$} \\
 \hline
 \textbf{Communication Topology} & fully connected & connected, directed & connected, undirected & {\color{black}even disconnected, directed} \\
 \hline
 \textbf{Suboptimality Guarantee} & {$1/2\,(\opt\,-\,\epsilon)$} & $1/2\,\opt$ & $1/2\,\opt$ &\!$1/2\,(\opt-{\sum_{i\in \calN} {\ourcurv_i}})$\!\\
 \hline
\end{tabular}
\caption{\textbf{\alg vs.~State of the Art.} The state of the art is divided into algorithms that optimize (i) in the continuous domain, employing a continuous representation of $f$~{\cite{calinescu2011maximizing}}, and (ii) in the discrete domain.  The continuous-domain algorithms need to compute the continuous representation's gradient via sampling; $M$ denotes the sample size ($M$ is $3$ in \cite{du2020jacobi}, $10$ in \cite{robey2021optimal}, and $1000$ in \cite{rezazadeh2021distributed}, in numerical evaluations with 10 or fewer agents).
}
\label{tab:comparison}
\end{table*}
}

We present \alg's computation, communication, and memory requirements. 
The requirements are in accordance to the paradigm in \Cref{tab:intro}, and are summarized in \Cref{tab:comparison}.  

We use the additional notation:
\begin{itemize}
    \item $\numl$ and $\el$ are the lengths of a message containing a real number or an action $s\in \calV_\calN$, respectively.
    \item $\scenario{diam}(\calG)$ is the \emph{diameter} of a network $\calG$, \ie the longest shortest path among any pair of nodes in $\calG$~{\cite{mesbahiBook}};
    \item $|\calX|$ is the cardinality of a discrete set $\calX$.
\end{itemize}

\begin{proposition}[Computation Requirements]\label{prop:computation}
Each agent  $i$ performs  $O(|\mathcal{V}_i||\mathcal{N}^-_i|)$ function evaluations during \alg. 
\end{proposition}
Each agent $i$ needs to re-evaluate the marginal gains of all $v\in\mathcal{V}_i$, every time an in-neighbor $j\in\mathcal{N}^-_i$ selects an action. Therefore, $i$ will perform $|\mathcal{N}^-_i||\mathcal{V}_i|$ evaluations in the worst case (and $|\mathcal{V}_i|$ evaluations in the best case).

\begin{proposition}[Communication Requirements]\label{prop:communication}
\alg's number of communication rounds is at most $2|\ouragent|-2$. 
\end{proposition}
Each iteration of \alg requires two communication rounds: one for marginal gains, and one for actions. 
Also, \alg requires $|\ouragent|-1$ iterations in the worst case (when only one agent selects an action at each iteration).  All in all, \alg requires at most $2|\ouragent|-2$ communication rounds.

\begin{proposition}[Memory Requirements]\label{prop:memory}
\alg's largest inter-agent message length is $\max{(\numl,\el)}$. 
\end{proposition}
Any inter-agent message in \alg contains either a marginal gain, or an action $s$.  Thus, the message's length is either $\numl$ or $\el$.  The total on-board memory requirements for each agent $i$ are $|\calN_i^-|\;\max{(\numl,\el)}$.

\begin{remark}[Near-Minimal Resource Requirements]
\label{rem:resource-awareness}
\alg has near-minimal computation, communication, and memory requirements, in accordance to \Cref{tab:intro}.  \emph{(i)  Computations per Agent:} the number of computations per agent is indeed proportional to the size of the agent's action set, and, in particular, is $O(|\mathcal{V}_i||\mathcal{N}^-_i|)$, \ie decreasing as $|\mathcal{N}^-_i|$ decreases.   The number of computations would have been minimal if instead it was $O(|\mathcal{V}_i|)$, since that is the cost for agent $i$ to compute its best action in $\calV_i$. \emph{(ii) Communication Rounds:} the number of communication rounds is indeed proportional to the number of agents, and, in particular, is at most $2|\calN|-2$.  The number is near-minimal, since, in the worst case of a line communication network, $|\calN|-1$ communication rounds are required for information to travel between the most distant agents. \emph{(iii) Memory per Message:} the length per message is indeed equal to the length of a real number or of an action.
\end{remark}
 
Besides, each agent $i$ can afford to run \alg, in accordance to \Cref{tab:intro}, by adjusting the size of its in- and out- neighborhoods $\calN_i^-$ and $\calN_i^+$.
\textit{E.g.}, by decreasing the size of $\calN_i^-$: (i) agent $i$'s computation effort decreases, since the effort is proportional to the size of $\calN_i^-$ (\Cref{prop:computation}); (ii) the per-round communication effort decreases, since the total communication rounds remain at most $2|\calN|-2$ 
(\Cref{prop:communication}) but the number of received messages per round decreases (\alg's lines 5--6); and (iii) the on-board memory-storage requirements decrease, since the inter-agent message length remains constant (\Cref{prop:memory}) but the number of received \mbox{messages per round decreases (\alg's lines 5--6).}

\begin{remark}[vs. State-of-the-Art Resource Requirements]\label{rem:comparison-resources}
\alg has comparable or superior computation, communication, and memory requirements, vs.~the state of the art.  
The comparison is summarized in \Cref{tab:comparison}. 

\myParagraph{Context and notation in \Cref{tab:comparison}} We divide the state of the art into algorithms that optimize (i) indirectly in the continuous domain, employing the continuous representation \emph{multi-linear extension}~\cite{calinescu2011maximizing} of the set function $f$~{\cite{robey2021optimal,du2020jacobi,rezazadeh2021distributed}}, and (ii) directly in the discrete domain~\cite{corah2018distributed, konda2021execution,liu2021distributed}.  The continuous-domain algorithms employ consensus-based techniques~{\cite{robey2021optimal,rezazadeh2021distributed}} or algorithmic game theory~{\cite{du2020jacobi}}, and %
require the computation of the multi-linear extension's gradient.  The computation is achieved via sampling; $M$ in Table~\ref{tab:comparison} denotes that sample size. $M$ is equal to  $3$ in \cite{du2020jacobi}, $10$ in \cite{robey2021optimal}, and $1000$ in \cite{rezazadeh2021distributed}, in numerical evaluations with 10 or fewer agents. 
The \emph{computations per agent} and \emph{communication rounds} reported for \cite{du2020jacobi} are based on the numerical evaluations therein, since a theoretical quantification is missing in \cite{du2020jacobi} and seems non-trivial to derive one as a function of $\calN$, $\epsilon$, or any other of the problem parameters. 
Further, all continuous-domain algorithms' resource requirements depend on additional problem-dependent parameters ({such as Lipschitz constants, the diameter of the domain set of the multi-linear extension, and a bound on the gradient of the multi-linear extension}), which here we make implicit via the big O, Omega, and Theta notation.  $\epsilon$ determines the approximation performance {of the respective algorithms.}

\myParagraph{Computations} Konda et al.~\cite{konda2021execution} rank best with $|\calV_i|$~computations per agent. \alg  ranks 2nd-best with $O(|\calN_i^-||\calV_i|)$.  
 The continuous-domain algorithms require a higher number of computations, proportional to $|\calN|^{2}$ or more.

\myParagraph{Communication} For undirected networks, Konda et al.~\cite{konda2021execution} and \alg rank best, requiring in the worst-case the same communication rounds; but \alg is also valid for directed networks. For appropriate $\epsilon$, the algorithm by Corah and Michael~\cite{corah2019distributed}, 
may require fewer communication rounds but in~\cite{corah2019distributed}, a pre-processing step with a \emph{fully connected network} is required. 
The remaining algorithms require a significantly higher number of communication rounds, proportional to $|\calN|^{2}$ or more.

\myParagraph{Memory} {\alg ranks best when $\scenario{length}_s\leq |\calV_\calN|\;\scenario{length}_{\#}$;  otherwise, it ranks after Robey et al.~\cite{robey2021optimal}~and Rezazadeh and Kia~\cite{rezazadeh2021distributed}, which then rank best (tie).} 
\end{remark}

%% file: 5-Guarantees.tex
\section{Approximation Guarantee of \alg: Centralization vs.~Decentralization Perspective}\label{sec:guarantee}

We present \alg's suboptimality bound (\Cref{th:additive}).  In accordance with \Cref{tab:intro}, the bound quantifies
the trade-off of \emph{centralization}, for global near-optimality, vs.~\emph{decentralization}, for near-minimal on-board resource requirements.

We introduce the notion of \emph{Centralization of Information among non-Neighbors} to quantify the bound. 

We also use the notation:
\begin{itemize}\setlength\itemsep{0.2em}
    \item $\mathcal{N}_i^c\triangleq\ouragent\setminus\{\mathcal{N}^-_i\cup\{i\}\}$ is the set of agents beyond the in-neighborhood of $i$ (see Fig.~\ref{fig:venn_diagram}), \ie $i$'s non-neighbors;
    \item $\mathcal{S}_{\calN_i^c}\triangleq \{s_j\}_{j\myin \calN_i^c}$ is the agents' actions in $\calN_i^c$;
    \item $\calS^\opt \myin\arg \,\max_{\single_i\myin\mathcal{V}_i, \, \forall\, i\myin \calN} \, f(\,\{s_i\}_{i\myin \calN}\,)$, \ie $\calS^\opt$ is an optimal solution to Problem~\ref{pr:main};
    \item $\oursol\triangleq\{\singlesol_i\}_{i\myin\ouragent}$ is \alg's output for all agents;
    \item $d(i,\calN_i^c)\triangleq\min_{j\myin\calN_i^c}\; d(i,j)$ for any agent $i\in\calN$ and distance metric $d$, \ie the distance between an agent $i$ and its non-neighborhood is the minimum distance between $i$ and a non-neighbor $j$.\\
\end{itemize}

\subsection{Centralization of Information:  A Novel Quantification}

We use the notion of \emph{Centralization Of Information among Non-neighbors} (\scenario{COIN}) to bound \alg's suboptimality. 

\definecolor{lightgray}{rgb}{0.83, 0.83, 0.83}
\begin{figure}[b]
\captionsetup{font=footnotesize}
\begin{center}
\vspace{-2mm}
\hspace{8mm}\begin{tikzpicture}
\filldraw[fill=lightgray] (-2,-2) rectangle (3,1);
\scope
\clip (-.25,-.5) circle (1);
\fill[white] (-.25,-.5) circle (1);
\endscope
\draw  (-.25,-.5) circle (1);
\draw  node at (-.25,-.5) {$\calN_i^-\; \cup \; \{i\}$};
\draw  node at (2,-.55) {$\calN_i^c$};
\draw  node at (3.3,0.85) {$\calN$};
\end{tikzpicture}
\end{center}
\caption{\textbf{Venn-diagram definition of the set $\calN_i^c$ for agent $i \in \calN$}. 
}\label{fig:venn_diagram}
\end{figure}

\input{5-figure-guarantees-image-covering}

\begin{definition}[Centralization Of Information among non-Neighbors (\scenario{COIN})]
\label{def:ourcurv}
Consider a communication network $\calG=(\ouragent,\mathcal{E})$, an agent $i\in \calN$, and a set function $f:2^{\calV_\calN}\mapsto$ $\mathbb{R}$. Then, agent $i$'s \emph{Centralization Of Information among non-Neighbors} is defined by 
\begin{equation}\label{eq:ourcurv}
  \ourcurv_i\triangleq {\max_{s_i\myin\mathcal{V}_i}\,\max_{s_j\myin\mathcal{V}_j,\,\forall\, j\myin \mathcal{N}_i^c}}\, \left[f(s_i) - f(s_i\,|\,\mathcal{S}_{\calN_i^c})\right].\!\!
\end{equation}
\end{definition}

If $f$ were entropy, then $\ourcurv_i$ looks like the mutual information between the information collected by agent $i$'s action $s_i$, and the information collected by agents $\calN_i^c$'s actions $\calS_{\calN_i^c}$, \ie the actions of agent $i$'s \emph{non-}neighbors. $\ourcurv_i$, in particular, is computed over the best such actions (hence the maximization with respect to $s_i$ and $\calS_{\calN_i^c}$ in \cref{eq:ourcurv}). %

$\ourcurv_i$ captures the \emph{centralization of information} within the context of a multi-agent network for information acquisition: $\ourcurv_i=0$ if and only if $s_i$'s information is independent from $\calS_{\calN_i^c}$, \ie if and only if the information is decentralized between agent $i$ and its non-neighbors $\calN_i^c$.

Computing $\ourcurv_i$ can be NP-hard~{\cite{Feige:1998:TLN:285055.285059}}, 
{or even impossible since agent $i$ may be unaware of its non-neighbors' actions}; 
but upper-bounding it can be easy.  In this paper, we focus on  upper bounds that depend on $d(i,\calN_i^c)$, quantifying how fast the information overlap between $i$ and its non-neighbors decreases the further away $\calN_i^c$ is from $i$; \ie how quickly information becomes decentralized beyond agent $i$'s neighborhood, the further away the non-neighborhood is.  For an \emph{image covering} task, we obtain such a bound next.

\begin{remark}[Distance-Based Upper Bounds for \scenario{COIN}: Image Covering Example]\label{rem:coin-bound}
Consider a toy image covering task
where each agent carries a camera with a \emph{round field of view of radius $r_s$} (Fig.~\ref{fig:image_covering}(a)).  Consider that each agent $i$ has fixed its in-neighbor $\calN_i^-$, \ie its communication range $r_i$ is just below the corresponding $d(i,\calN_i^c)$, that is, $$r_i=d(i,\calN_i^c)-\delta$$ for some, possibly arbitrarily small, $\delta>0$; \eg in Fig.~\ref{fig:image_covering}(b) the non-neighbor on the left of agent $i$ is just outside the boundary of $i$'s communication range.
Then, $\ourcurv_i$ is equal to the overlap of the field of views of agent $i$ and its non-neighbors, assuming, for simplicity, that the bound remains the same across two consecutive moves.  Since the number of agent $i$'s non-neighbors may be unknown, 
an upper bound to $\ourcurv_i$ is the gray ring area in Fig.~\ref{fig:image_covering}(b), obtained assuming an infinite amount of non-neighbors around agent $i$, located just outside the boundary of $i$'s communication range.  That~is, 
\begin{equation}\label{eq:communication_range}
\ourcurv_i \leq \max(0, \;\pi [r_s^2-(r_i-r_s)^2]). %
\end{equation} 
The upper bound in \cref{eq:communication_range} as a function of the communication range $r_i$ is shown in Fig.~\ref{fig:image_covering}(c). As expected, it tends to zero for increasing $r_i$, equivalently, for increasing $d(i,\calN_i^c)$.  Particularly, when $d(i,\calN_i^c)>2r_s$, the field of views of robot $i$ and of each of the non-neighboring robots $\calN_i^c$ are non-overlapping, and, thus, $\ourcurv_i=0$.
\end{remark}

In \Cref{subsec:approx-guarantee}, we show that \alg enables each agent~$i$ to choose its in-neighborhood $\calN_i^-$ (equivalently, its non-neighborhood $\calN_i^c$) to balance both its on-board resources and its contribution to a near-optimal approximation performance, as the latter is captured by $\ourcurv_i$.

\begin{remark}[Relation to~Pairwise Redundancy{~\cite{corah2018distributed}}]\label{rem:pairwise}
$\ourcurv$'s definition generalizes the notion of \emph{pairwise redundancy} $w_{ij}$ between two agents $i$ and $j$, introduced by Corah and Michael~\cite{corah2018distributed}.  The notion was introduced in the context of parallelizing the execution of the sequential greedy~\cite{fisher1978analysis}, by ignoring the edges between pairs of agents in an a priori \emph{fully connected network}.  The comparison of the achieved parallelized greedy in~\cite{corah2018distributed} and \alg is found in \Cref{tab:comparison}.   Besides, $w_{ij}$ captures the mutual information between the two agents $i$ and $j$, defined as $w_{ij}\triangleq \max_{s_i\myin\mathcal{V}_i}\, \max_{s_j\myin\mathcal{V}_j,\,j\myin\mathcal{N}_i}\,[f(s_i)-f(s_i\,|\,s_j)]$; whereas $\ourcurv_i$ captures the mutual information between an agent $i$ and \emph{all} its \emph{non}-neighbors, capturing directly the decentralization of information across the network.  
\end{remark}

\subsection{Approximation Guarantee of \alg}\label{subsec:approx-guarantee}

\begin{theorem}[Approximation Performance of \alg]\label{th:additive}
\alg selects $\oursol$ such that $\singlesol_i\myin\mathcal{V}_i, \, \forall\, i\myin \calN$, and 
{
\beq\label{eq:approx}
   f(\oursol)\,\geq \, \frac{1}{2}\,\left[f(\mathcal{S}^\opt) \, -\, \sum_{i\myin \calN}\ourcurv_i\right].
\eeq}
\end{theorem}

\paragraph*{Proof of \Cref{th:additive}} We index each agent in $\calN$ per its selecting order in \alg, \ie  agent $i \in \ouragent\triangleq\{1,\dots,|\ouragent|\}$ is the $i$-th agent to select an action during the execution of \alg. If two agents select actions simultaneously, then we index them randomly.  We use also the notation:
\begin{itemize}\setlength\itemsep{0.2em}
    \item  $[i]$ is the set of agents $\{1, \dots, i\}$;
    \item $\calS^{\scenario{RAG}}_\calX \triangleq \{s^\alg_i\}_{i\myin \calX}$ for any $\calX\subseteq \calN$, \ie $\calS^{\scenario{RAG}}_\calX$ is the set of actions selected by the agents in $\calX$.
\end{itemize}

The first 3 steps in the proof (\cref{aux1:1,aux1:2,aux1:3}) are inspired by~\cite{fisher1978analysis}; the 5-th and 6-th steps (\cref{aux1:5,aux1:7}) are inspired by~\cite{corah2018distributed}.  All steps involve novel quantities and result in a novel suboptimality bound for the novel algorithm \alg. 
{%
\begin{align}
    &f(\mathcal{S}^\opt)\nonumber\\
    &\leq f(\mathcal{S}^\opt,\,\oursol) \label{aux1:1} \\
    &= f(\oursol) + \sum_{i \myin \calN} \,f(s_i^\opt\,|\,\oursol,\,\mathcal{S}^\opt_{[i-1]}) \label{aux1:2} \\
    &\leq f(\oursol) + \sum_{i \myin \calN}\, f(s_i^\opt\,|\,\oursol_{\mathcal{N}_i^-\cap[i-1]}) \label{aux1:3} \\
    &\leq f(\oursol) + \sum_{i \myin \calN}\, f(s_i^\alg\,|\,\oursol_{\mathcal{N}_i^-\cap[i-1]}) \label{aux1:4} \\
    &= 2f(\oursol) + \sum_{i \myin \calN} \Big[f(s_i^\alg\,|\,\oursol_{\mathcal{N}_i^-\cap[i-1]})\nonumber\\
    &\quad\quad\quad\quad\quad\quad\quad\quad\quad\quad\quad\quad\quad\;\;-f(s_i^\alg\,|\,\oursol_{[i-1]})\Big]  \label{aux1:5}\\
    &= 2f(\oursol) + \sum_{i \myin \calN} \Big[f(s_i^\alg\,|\,\oursol_{\mathcal{N}_i^-\cap[i-1]}) \nonumber\\
    &\quad\quad\quad\quad\quad\quad\quad- f(s_i^\alg\,|\,\oursol_{\mathcal{N}_i^-\cap[i-1]},\,\oursol_{[i-1]\setminus\mathcal{N}_i^-})\Big]\label{aux1:6} \\
    &\leq 2f(\oursol) + \sum_{i \myin \calN} \left[f(s_i^\alg) - f(s_i^\alg\,|\,\oursol_{[i-1]\setminus\mathcal{N}_i^-})\right] \label{aux1:7} \\
    &\leq 2f(\oursol) + \sum_{i \myin \calN} \left[f(s_i^\alg) - f(s_i^\alg\,|\,\oursol_{\calN_i^c})\right] \label{aux1:8} \\
    &\leq 2f(\oursol) + \sum_{i \myin \calN}\, \ourcurv_i,\label{aux1:9} 
\end{align}}where \cref{aux1:1} holds true due to the monotonicity of $f$; \cref{aux1:2,aux1:5} are proved by telescoping the sums; \cref{aux1:3,aux1:8} hold true due to the submodularity of $f$; \cref{aux1:4} holds true since \alg selects $\singlesol_i$ greedily; \cref{aux1:6} holds true since $\oursol_{[i-1]}=\oursol_{\mathcal{N}_i^-\cap[i-1]}\cup\oursol_{[i-1]\setminus\mathcal{N}_i^-}$%
; \cref{aux1:7} holds true due to the 2nd-order submodularity; and \cref{aux1:9} holds true due to $\ourcurv_i$'s definition (Definition~\ref{def:ourcurv}). \qed

\begin{remark}[\textit{Centralization} vs.~\textit{Decentralization}]\label{rem:trade-off}
\alg's suboptimality bound in \Cref{th:additive} captures the trade-off of \emph{centralization}, for global near-optimality, vs.~\emph{decentralization}, for near-minimal on-board resource requirements:
\begin{itemize}
    \item \emph{Near-optimality} requires \emph{large} $\calN_i^-$, \ie \emph{centralization}: the larger $\calN_i^-$ is, the larger $d(i,\calN_i^c)$ is. Thus, the smaller is the information overlap between $i$ and its non-neighbors $\calN_i^c$, \ie the smaller $\ourcurv_i$ is, resulting in an increased near-optimality for \alg.
    \item \emph{Minimal on-board resource requirements} requires instead a \emph{small} $\calN_i^-$, \ie \emph{decentralization}:  the smaller $\calN_i^-$ is, the less the computation per agent (\Cref{prop:computation}), as well as, the less the per-communication-round communication and memory-storage effort, since the number of received messages per round decreases. 
\end{itemize} 
\end{remark}

\alg covers the spectrum from fully centralized  ---all agents communicate with all others--- to fully decentralized ---all agents communicate with none. \alg enjoys the suboptimality guarantee in \Cref{th:additive} throughout the spectrum, capturing the trade-off of centralization vs.~decentralization. When fully centralized, \alg matches the  $1/2$ suboptimality bound of the classical greedy~\cite{fisher1978analysis}, since then $\calN_i^c=\emptyset$ for all $i$, \ie $\ourcurv_i=0$.  For a centralized algorithm, the best possible bound is $1-1/e\simeq .63$~\cite{sviridenko2017optimal}.

\begin{remark}[vs. State-of-the-Art Approximation Guarantees]\label{rem:comparison-approximation}
\alg is the first  algorithm to quantify the trade-off
of \emph{centralization \emph{vs.}~decentralization}, per \Cref{tab:intro}.  \alg enables each agent to independently decide the size of its in-neighborhood to balance \emph{near-optimality} ---which requires larger in-neighborhoods, \ie \emph{centralization}--- and \emph{on-board resources} ---which requires smaller in-neighborhoods, \ie \emph{decentralization}.  
Instead, the state of the art tunes near-optimality via a globally known hyper-parameter $\epsilon$, without accounting for the balance of smaller vs.~larger neighborhoods at the agent level, and independently for each agent.

{Moreover, the continuous methods \cite{robey2021optimal,du2020jacobi,rezazadeh2021distributed} achieve the suboptimality bound in probability,
and the achieved value is in expectation.} Instead, \alg's  bound is deterministic, involving the exact value of the selected actions.
\end{remark}

In the \emph{image covering} scenario of \Cref{rem:coin-bound}, \alg's suboptimality guarantee gradually degrades with increased decentralization: as the agent $i$'s {communication range} decreases,
$d(i,\calN_i^c)$ becomes smaller, and $\ourcurv_i$ increases (Fig.~\ref{fig:image_covering}(b)), causing \alg's suboptimality guarantee to decrease.

%% file: 5-figure-guarantees-image-covering.tex
\begin{figure*}[t]
    \captionsetup{font=footnotesize}
	\begin{center}
	\hspace{-9cm}\begin{minipage}{\columnwidth}
	\hspace{-2mm}		
	    \begin{minipage}{0.6\columnwidth}%

            \centering%
            \includegraphics[width=.9\columnwidth]{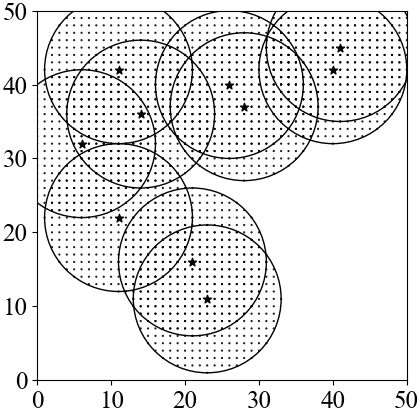} \\
            \caption*{(a) \textbf{Image covering scenario in a $50\text{ points}$ $ \times 50$ points map with $10$ agents.} The stars are agents' locations, the circles are agents' sensing ranges, and the dots covered points.
            }
		\end{minipage}~~~~~
		\begin{minipage}{0.75\columnwidth}%
	    \vspace{0cm}
        \centering%
        \includegraphics[width=.95\columnwidth]{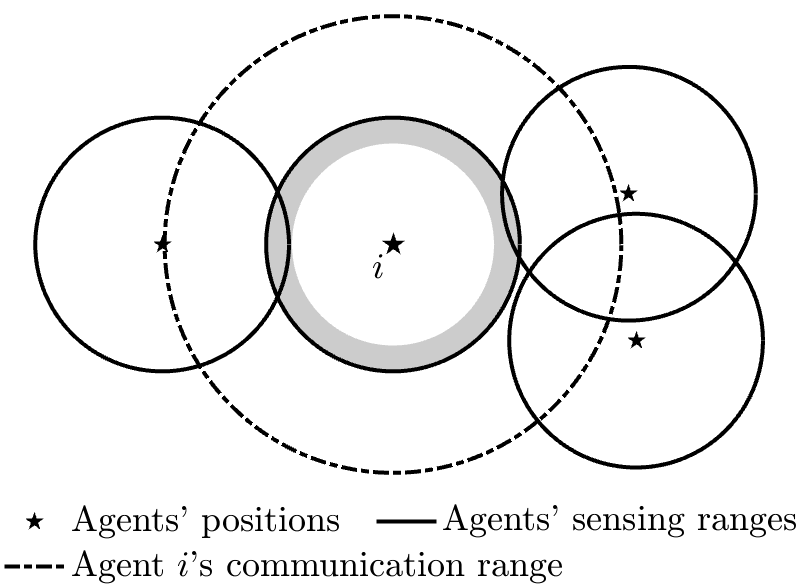} \\ \caption*{(b) \textbf{Agent $i$ and its non-neighbors}. Non-neighbors are the agents beyond agent $i$'s communication range. The sensing ranges of each agent are also depicted, defining the area of their field of view. 
        }
		\end{minipage}~~~~~
		\begin{minipage}{0.6\columnwidth}%
			    \vspace{-.3cm}
            \centering%
            \hspace*{-2mm}\includegraphics[width=1\columnwidth]{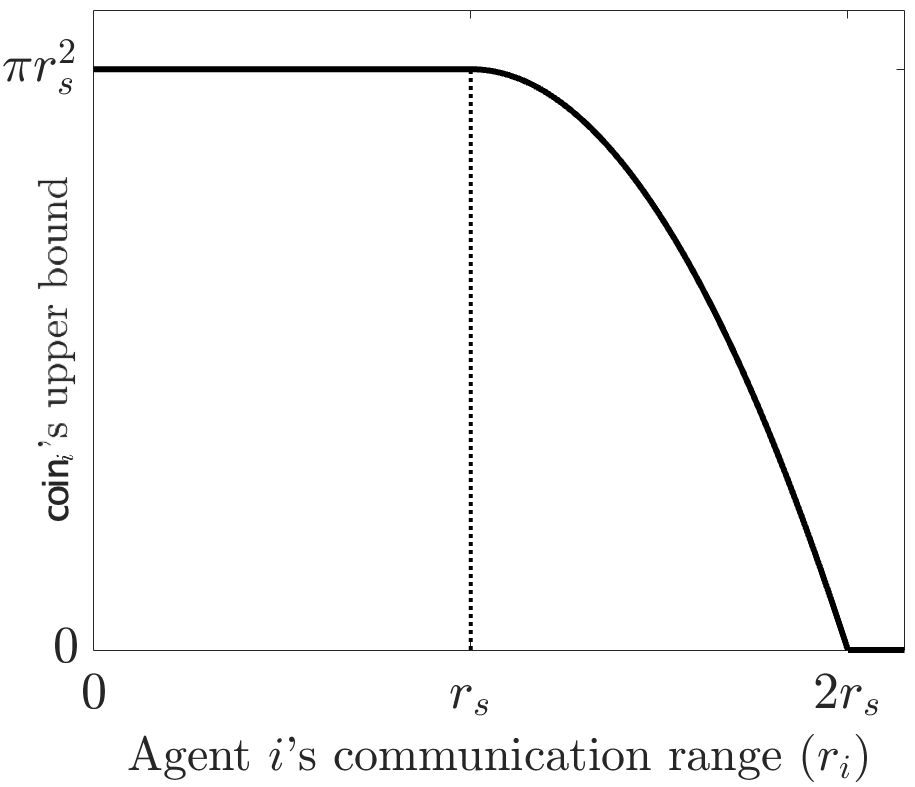} \\
             \vspace{0.4mm}
            \caption*{(c) \textbf{$\ourcurv_i$'s upper bound for increasing agent $i$'s communication range}.  $r_s$ is the agents' sensing range.%
            }
		\end{minipage}
	\end{minipage}
	\caption{\textbf{Image Covering Setup}. 
	(a) A scenario;  (b) Agent $i$ and its non-neighbors; (c)  $\ourcurv_i$'s upper bound for increasing agent $i$'s communication range.
	}\label{fig:image_covering}
	\end{center}
\end{figure*}

%% file: 6-Experiments.tex
\section{Evaluation in Image Covering with Robots}\label{sec:experiments}

We evaluate \alg in simulated scenarios of image covering with mobile robots (Fig.~\ref{fig:image_covering}). We first compare \alg with the state of the art (\Cref{subsec:vs-covering}; see \Cref{tab:sim-comparison}).  Then, we evaluate the trade-off of centralization vs.~decentralization with respect to \alg's performance (\Cref{subsec:trade-off-covering}; see Fig.~\ref{fig:different_comm_range}).

We performed all simulations in Python 3.9.7, on a MacBook Pro with the Apple M1 Max chip and a 32 GB RAM. 

Our code is open-sourced here: \href{https://github.com/UM-iRaL/Resource-Aware-Coordination-AirSim}{\color{red}https://github.com/UM-iRaL/Resource-Aware-Coordination-AirSim}.

\subsection{\alg vs.~State of the Art}\label{subsec:vs-covering}

We compare \alg with the state of the art in simulated scenarios of \textit{image covering} (Fig.~\ref{fig:image_covering}). To enable the comparison, %
we set up undirected and connected communication networks (\alg is valid for directed and even disconnected networks but the state-of-the-art methods are \textit{not}).  \alg demonstrates superior or comparable performance (\Cref{tab:sim-comparison}).

\myParagraph{Simulation Scenario} 
{We consider 50 instances of the setup in Fig.~\ref{fig:image_covering}(a). Without loss of generality, each agent has a communication range of $15$, a sensing radius of $10$, and the action set \!\{``forward'', ``backward'', ``left'', ``right''\}\! by 1 point. {The agents seek to maximize the number of covered~points.}}

\myParagraph{Compared Algorithms}
{We compare \alg with the methods by Robey et al.~\cite{robey2021optimal} and Konda et al.~\cite{konda2021execution} since, among the state of the art, they achieve top performance for at least one resource requirement and/or for their suboptimality guarantee (\Cref{tab:comparison}) ---the method by Corah and Michael \cite{corah2018distributed} may achieve  fewer communication rounds, yet it requires a fully connected network.
To ensure the method in~\cite{robey2021optimal} achieves a sufficient number of covered points, we set the sample size $M=10$, as is also set in~\cite{robey2021optimal}, and the number of communication rounds $T=100$.
}

\myParagraph{Results} The results are reported in \Cref{tab:sim-comparison}, and mirror the theoretical comparison in \Cref{tab:comparison}.
\alg demonstrates superior or comparable performance, requiring, (i) $5$ orders of magnitude less computation time vs.~the state-of-the-art consensus algorithm in~\cite{robey2021optimal}, and comparable computation time vs.~the state-of-the-art greedy algorithm in~\cite{konda2021execution}, (ii)~$1$ order of magnitude fewer communication rounds vs.~the consensus algorithm in~\cite{robey2021optimal},  and comparable communication rounds vs.~the greedy algorithm in~\cite{konda2021execution}, and (iii) the least memory (\eg $40\%$ less than the consensus algorithm).  Still, \alg achieves the best approximation performance.

{\footnotesize
\begin{table}[t]
\captionsetup{font=footnotesize}
\centering
\renewcommand{\arraystretch}{1.5}
\vspace*{2.5mm}
\begin{tabular}{ |l||c|c|c| }
 \hline
 \!\!\textbf{Method} & \!\!\begin{tabular}{@{}c@{}}
  {Robey et al.}\vspace{-1mm}\\ \cite{robey2021optimal}   \end{tabular}\!\!& \!\!\begin{tabular}{@{}c@{}}
  {Konda et al.}\vspace{-1mm}\\ \cite{konda2021execution}
  \end{tabular}\!\!&\!\!\begin{tabular}{@{}c@{}}
  {\alg}\vspace{-1mm}\\ (this paper)
  \end{tabular}\!\!\\
 \hline 
 \hline
 \!\!\textbf{Total Computation Time (s)}\!\!\!
 & 1434.11 &	\textbf{0.02}	& 0.05\\
 \hline
 \!\!\textbf{Communication Rounds} & 100 & 13.44 & \textbf{7.76} \\
 \hline
 \!\!\textbf{Peak Total Memory (MB)} & 290.43 &	181.07&	\textbf{167.07}\\
 \hline
 \!\!\textbf{Total Covered Points} &1773.54& 	1745.18	& \textbf{1816.4}\\
 \hline
\end{tabular}
\renewcommand{\arraystretch}{1}
\caption{\textbf{\alg vs.~State of the Art.} Averaged performance over $50$ image covering instances  involving 10 robots in a $50\text{ points} \times 50$ points map. 
}
\label{tab:sim-comparison}
\end{table}
}

\subsection{The Trade-Off of Centralization vs.~Decentralization}\label{subsec:trade-off-covering}

We demonstrate the trade-off of \textit{centralization}
vs.~\textit{decentralization}, 
with respect to \alg's performance.

\myParagraph{Simulation Scenario} {We consider the same setup as in ~\Cref{subsec:vs-covering}, yet with the communication range increasing from $1$ to $50$. That is,} the communication network starts from being fully disconnected (fully decentralized) and becomes fully connected (fully centralized). %
The communication range is assumed the same for all robots, for simplicity.

\begin{figure}[t]
    \captionsetup{font=footnotesize}
    \centering
    \includegraphics[width=.98\columnwidth]{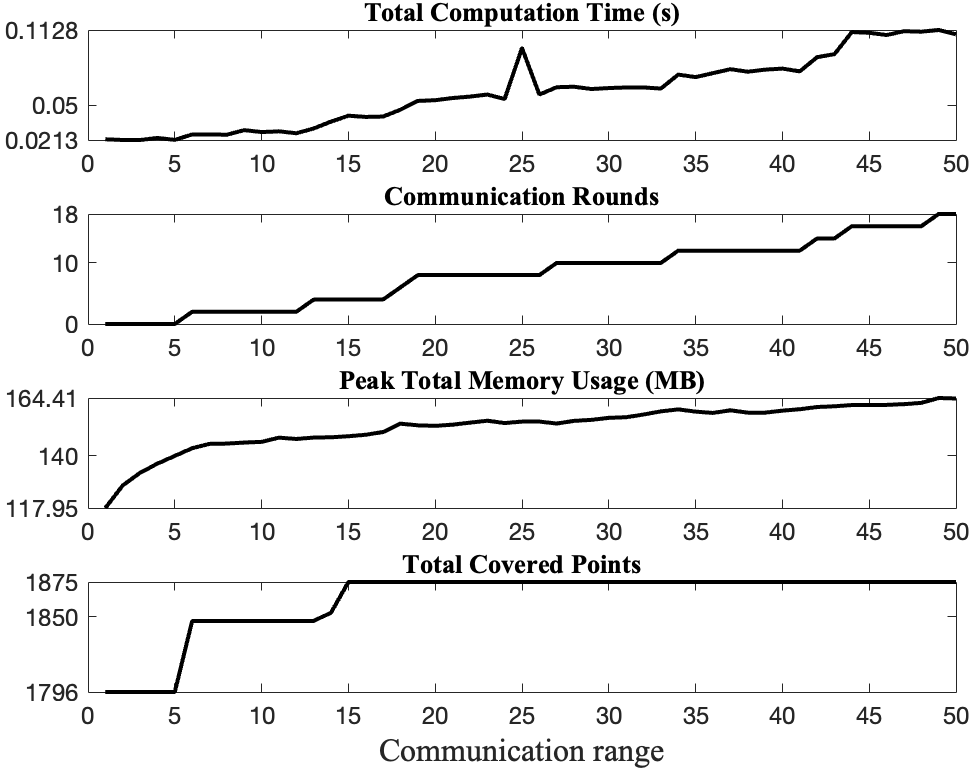}
        \vspace{-1mm}
    \caption{\textbf{Centralization vs.~Decentralization: Resource requirements and coverage performance of \alg for increasing communication range}, in an image covering scenario with $10$ robots in a $50\text{ points} \times 50$ points map. 
    } 
    \label{fig:different_comm_range}
\end{figure}

\myParagraph{Results}
The results are reported in Fig.~\ref{fig:different_comm_range}. When a higher communication range results in more in-neighbors, then (i) 
each agent 
executes more iterations of \alg before selecting an action, resulting in increased (i-a) 
computation time and (i-b) communication rounds ($1$ iteration of \alg corresponds to $2$ communication rounds; see \alg's lines 5 and 11), and (ii) each agent needs more on-board memory for information storage and processing. In contrast, with more in-neighbors, each agent coordinates more centrally, and, thus, the total covered points increase  (for each agent $i$, $\ourcurv_i$ becomes smaller till it vanishes, and the theoretical suboptimality bound  \mbox{increases to $1/2$). All in all,} \textit{Fig.~\ref{fig:different_comm_range} captures the trade-off of \emph{centralization}, for global near-optimality, vs.~\emph{decentralization}, for near-minimal on-board resource requirements.  For increasing communication range, the required on-board resources increase, but also the total covered points increase. To balance the trade-off, the communication range may be set to~$6$.}

%% file: 7-Conclusion.tex
\section{Conclusion} \label{sec:con}
\myParagraph{Summary} We made the first step to enable a resource-aware distributed intelligence among heterogeneous agents.  
We are motivated by complex tasks taking the form of \Cref{pr:main}, such as image covering. 
We introduced a \textit{resource-aware optimization paradigm} (\Cref{tab:intro}),
and presented \alg, the first resource-aware algorithm.
\alg is the first algorithm to quantify the trade-off of \emph{centralization}, for global near-optimality, vs.~\emph{decentralization}, for near-minimal on-board resource requirements.
To capture the trade-off, we introduced the notion of \emph{Centralization of Information among non-Neighbors} (\scenario{COIN}).  
We validated \alg in simulated scenarios of image covering, demonstrating its superiority. 

\myParagraph{Future Work} \alg assumes synchronous communication.  Besides, the communication topology has to be fixed and failure-free across communication rounds (\Cref{ass:fixed}).  
Our future work will enable \alg beyond the above limitations.
We will also consider multi-hop communication.  Correspondingly, we will quantify the trade-off  of~near-optimality vs.~resource-awareness
based on the depth of the multi-hop communication.  We will also extend our results to any submodular function (instead of 2nd-order submodular). 

Further, we will leverage our prior work~\cite{tzoumas2022robust} to extend \alg to the attack-robust case, against robot removals.

%% file: references.bib
@article{du2020jacobi,
  title={Jacobi-Style Iteration for Distributed Submodular Maximization},
  author={Du, Bin and Qian, Kun and Claudel, Christian and Sun, Dengfeng},
  journal={arXiv:2010.14082},
  year={2020}
}

@article{mcguire2019minimal,
  title={Minimal navigation solution for a swarm of tiny flying robots to explore an unknown environment},
  author={McGuire, KN and De Wagter, Christophe and Tuyls, Karl and Kappen, HJ and de Croon, Guido CHE},
  journal={Science Robotics},
  volume={4},
  number={35},
  pages={9710},
  year={2019},
  publisher={American Association for the Advancement of Science}
}

@article{tzoumas2022robust,
  title={Robust and adaptive sequential submodular optimization},
  author={Tzoumas, Vasileios and Jadbabaie, Ali and Pappas, George J},
  journal={IEEE Transactions on Automatic Control},
    volume={67},
  number={1},
    pages={89-104},
  year={2022},
}

@article{gharesifard2017distributed,
  title={Distributed submodular maximization with limited information},
  author={Gharesifard, Bahman and Smith, Stephen L},
  journal={IEEE Transactions on Control of Network Systems (TCNS)},
  volume={5},
  number={4},
  pages={1635--1645},
  year={2017},
}

@article{mirzasoleiman2016distributed,
  title={Distributed submodular maximization},
  author={Mirzasoleiman, Baharan and Karbasi, Amin and Sarkar, Rik and Krause, Andreas},
  journal={The Journal of Machine Learning Research (JMLR)},
  volume={17},
  number={1},
  pages={8330--8373},
  year={2016},
  publisher={JMLR. org}
}

@article{rezazadeh2021distributed,
  title={Distributed Strategy Selection: {A} Submodular Set Function Maximization Approach},
  author={Rezazadeh, Navid and Kia, Solmaz S},
  journal={arXiv preprint:2107.14371},
  year={2021}
}

@inproceedings{robey2021optimal,
  title={Optimal algorithms for submodular maximization with distributed constraints},
  author={Robey, Alexander and Adibi, Arman and Schlotfeldt, Brent and Hassani, Hamed and Pappas, George J},
  booktitle={Learning for Dynamics \& Control},
  pages={150--162},
  year={2021},
}

@article{liu2021distributed,
  title={Distributed resilient submodular action selection in adversarial environments},
  author={Liu, Jun and Zhou, Lifeng and Tokekar, Pratap and Williams, Ryan K},
  journal={IEEE Robotics and Automation Letters},
  volume={6},
  number={3},
  pages={5832--5839},
  year={2021},
  publisher={IEEE}
}

@article{corah2019distributed,
  title={Distributed matroid-constrained submodular maximization for multi-robot exploration: Theory and practice},
  author={Corah, Micah and Michael, Nathan},
  journal={Autonomous Robots (AURO)},
  volume={43},
  number={2},
  pages={485--501},
  year={2019},
}

@inproceedings{corah2018distributed,
  title={Distributed submodular maximization on partition matroids for planning on large sensor networks},
  author={Corah, Micah and Michael, Nathan},
  booktitle={IEEE Conference on Decision and Control (CDC)},
  pages={6792--6799},
   year={2018}
}

@article{konda2021execution,
  title={Execution Order Matters in Greedy Algorithms with Limited Information},
  author={Konda, Rohit and Grimsman, David and Marden, Jason},
  journal={arXiv preprint:2111.09154},
  year={2021}
}

@article{sviridenko2017optimal,
  title={Optimal approximation for submodular and supermodular optimization with bounded curvature},
  author={Sviridenko, Maxim and Vondr{\'a}k, Jan and Ward, Justin},
  journal={Math.~of Operations Research},
  volume={42},
  number={4},
  pages={1197--1218},
  year={2017},
}

@incollection{fisher1978analysis,
  title={An analysis of approximations for maximizing submodular set functions -- {II}},
  author={Fisher, Marshall L and Nemhauser, George L and Wolsey, Laurence A},
  booktitle={Polyhedral combinatorics},
  pages={73--87},
  year={1978},
}

@article{chen2022higher,
  title={Higher order monotonicity and submodularity of influence in social networks: from local to global},
  author={Chen, Wei and Li, Qiang and Shan, Xiaohan and Sun, Xiaoming and Zhang, Jialin},
  journal={Information and Computation},
  pages={104864},
  year={2022},
  publisher={Elsevier}
}

@article{wang2015accelerated,
  title={An accelerated continuous greedy algorithm for maximizing strong submodular functions},
  author={Wang, Zengfu and Moran, Bill and Wang, Xuezhi and Pan, Quan},
  journal={J.~of Combinatorial Optimization},
  volume={30},
  number={4},
  pages={1107--1124},
  year={2015},
  publisher={Springer}
}

@article{ramalingam2017efficient,
  title={Efficient minimization of higher order submodular functions using monotonic boolean functions},
  author={Ramalingam, Srikumar and Russell, Chris and Ladick{\'y}, LU and Torr, Philip HS},
  journal={Discrete Applied Math.},
  volume={220},
  pages={1--19},
  year={2017},
  publisher={Elsevier}
}

@article{downie2022submodular,
  title={Submodular Maximization with Limited Function Access},
  author={Downie, Andrew and Gharesifard, Bahman and Smith, Stephen L},
  journal={arXiv preprint:2201.00724},
  year={2022}
}

@incollection{crama1989characterization,
  title={A characterization of a cone of pseudo-Boolean functions via supermodularity-type inequalities},
  author={Crama, Yves and Hammer, Peter L and Holzman, Ron},
  booktitle={Quantitative Methoden in den Wirtschaftswissenschaften},
  pages={53--55},
  year={1989},
  publisher={Springer}
}

@article{foldes2005submodularity,
  title={Submodularity, supermodularity, and higher-order monotonicities of pseudo-Boolean functions},
  author={Foldes, Stephan and Hammer, Peter L},
  journal={Mathematics of Operations Research},
  volume={30},
  number={2},
  pages={453--461},
  year={2005},
  publisher={INFORMS}
}

@book{camacho2013model,
  title={Model {P}redictive {C}ontrol},
  author={Camacho, Eduardo F and Alba, Carlos Bordons},
  year={2013},
  publisher={Springer}
}

@article{krause08efficient,
	Author = {Andreas Krause and Jure Leskovec and Carlos Guestrin and Jeanne VanBriesen and Christos Faloutsos},
	Journal = {Journal of Water Resources Planning and Management},
	Number = {6},
	Pages = {516-526},
	Title = {Efficient Sensor Placement Optimization for Securing Large Water Distribution Networks},
	Volume = {134},
	Year = {2008}
}

@inproceedings{roberts2017submodular,
  title={Submodular trajectory optimization for aerial 3{D} scanning},
  author={Roberts, Mike and Dey, Debadeepta and Truong, Anh and Sinha, Sudipta and Shah, Shital and Kapoor, Ashish and Hanrahan, Pat and Joshi, Neel},
  booktitle={Proceedings of the IEEE International Conference on Computer Vision (ICCV)},
  pages={5324--5333},
  year={2017}
}

@inproceedings{mokhtari2018decentralized,
  title={Decentralized submodular maximization: Bridging discrete and continuous settings},
  author={Mokhtari, Aryan and Hassani, Hamed and Karbasi, Amin},
  booktitle={International Conference on Machine Learning},
  pages={3616--3625},
  year={2018},
}

@inproceedings{atanasov2015decentralized,
  title={Decentralized active information acquisition: {T}heory and application to multi-robot {SLAM}},
  author={Atanasov, Nikolay and Le Ny, Jerome and Daniilidis, Kostas and Pappas, George J},
  booktitle={IEEE International Conference on Robotics and Automation (ICRA)},
  pages={4775--4782},
  year={2015},
}

@inproceedings{corah2021scalable,
  title={Scalable Distributed Planning for Multi-Robot, Multi-Target Tracking},
  author={Corah, Micah and Michael, Nathan},
  booktitle={IEEE/RSJ International Conference on Intelligent Robots and Systems (IROS)},
  pages={437--444},
  year={2021}
}

@article{kumar2004robot,
  title={Robot and sensor networks for first responders},
  author={Kumar, Vijay and Rus, Daniela and Singh, Sanjiv},
  journal={IEEE Pervasive computing},
  volume={3},
  number={4},
  pages={24--33},
  year={2004},
  publisher={IEEE}
}

@article{calinescu2011maximizing,
  title={Maximizing a monotone submodular function subject to a matroid constraint},
  author={Calinescu, Gruia and Chekuri, Chandra and P{\'a}l, Martin and Vondr{\'a}k, Jan},
  journal={SIAM Journal on Computing},
  volume={40},
  number={6},
  pages={1740--1766},
  year={2011},
  publisher={SIAM}
}

@article{Feige:1998:TLN:285055.285059,
 author = {Feige, Uriel},
 title = {A Threshold of $ln(n)$ for Approximating Set Cover},
 journal = {Journal of the ACM},
 volume = {45},
 number = {4},
 year = {1998},
 pages = {634--652},
}

@article{grimsman2018impact,
  title={The impact of information in distributed submodular maximization},
  author={Grimsman, David and Ali, Mohd Shabbir and Hespanha, Joao P and Marden, Jason R},
  journal={IEEE Transactions on Control of Network Systems (TCNS)},
  volume={6},
  number={4},
  pages={1334--1343},
  year={2018},
}

@BOOK{mesbahiBook,
    author = {Mehran Mesbahi and Magnus Egerstedt},
    publisher = {Princeton University Press},
    title = {Graph Theoretic Methods in Multiagent Networks},
    year = {2010},
}
